# New Rigid-foldable Developable Quadrilateral Creased Papers

Zeyuan He[a]*, Simon D. Guest[a]

[a]Department of Engineering, University of Cambridge
Trumpington Street, Cambridge CB2 1PZ, UK
*zh299@cam.ac.uk

**Abstract**

This article extends the range of 1-DOF rigid-foldable developable quadrilateral creased papers. In a previous article, we put forward a sufficient and necessary condition for a quadrilateral creased paper to be rigid-foldable, and introduce a special sufficient condition that is convenient for practical use, generating quadrilateral creased paper by stitching basic units. In this article we develop new flat-foldable units and show how these can be used to construct a series of more complex rigid-foldable developable quadrilateral creased papers.

**Keywords**: rigid origami, rigid-foldability, isometric transformation, degree-4

## 1. Introduction

Rigid origami is a branch of origami in which it is possible for a paper to fold and unfold by bending along creases without deformation of the panels. It has been applied to many engineering areas with immensely different length scales, such as DNA origami, catalysts, biomedical devices, robotics, stretchable electronics, stacked meta-materials, architecture and aerospace elements (Callens and Zadpoor [7], Peraza-Hernandez et al. [8]).

Currently, most applications are based on relatively simple crease patterns, which inspires us to find further rigid-foldable creased papers. The sufficient and necessary condition for a developable quadrilateral creased paper that is flat-foldable to be rigid-foldable has been given in Tachi [3]. In our previous work [2], we have given a wider sufficient and necessary condition for a general quadrilateral creased paper to be rigid-foldable. We also give a useful sufficient condition by "stitching" together "basic units". In this article, we extend the set of possible units that can be used to generate rigid-foldable quadrilateral creased paper by describing new flat-foldable units. Note that we have provided the underlying definitions for the technical terms in [1]. Before further discussion we will first clarify the notion of a developable quadrilateral creased paper.

**Definition 1** A *quadrilateral* creased paper $Q$ satisfies:

1) $Q$ only contains degree-4 single-vertex creased papers.
2) All the inner panels are quadrilateral.

We say a *degree-4 single-vertex creased paper* is the union of four panels around an inner vertex and the crease pattern on the boundaries of them. A *developable* creased paper has the sum of sector angles of every degree-4 single-vertex creased paper equal to $2\pi$, and a rigidly folded state where all folding angles are 0, called the trivial rigidly folded state.

In section 2, the configuration space of a developable degree-4 single-vertex creased paper is proposed. This is the basic geometry used here to construct a rigid-foldable creased paper. We define a "unit" as the union of two degree-4 single-vertex creased papers sharing two panels and has special relations among its folding angles. By presenting a sufficient and necessary condition for rigid-foldability in section 3, we derive the relation among the sector angles of a flat-foldable unit in section 4. We can then stitch units together to form a series of rigid-foldable creased papers. Finally, we give some examples.

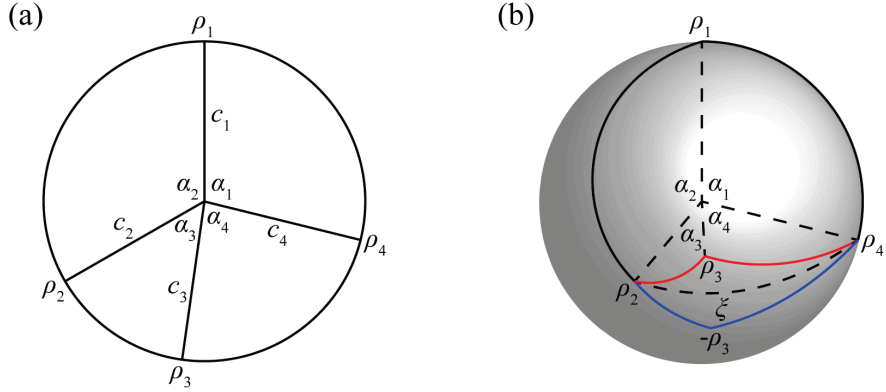

Figure 1: (a) A degree-4 single-vertex creased paper. We label these sector angles and corresponding inner creases counter-clockwise. (b) Two non-trivial rigidly folded states with the outside edge of the single-vertex creased paper drawn on a sphere as arcs of great circles, assuming the panel corresponding to $\alpha_1$ is fixed when we change the magnitude of $\rho_1$. Generically, there are two sets of folding angles $\rho_2, \rho_3, \rho_4$, which are symmetric to $\xi$ and colored by red and blue.

In this article we only consider a 1-connected paper. Paper with holes is discussed in [2].

## 2. Developable Degree-4 Single-vertex Creased Paper

A developable degree-4 single-vertex creased paper has four sector angles $\alpha_1, \alpha_2, \alpha_3, \alpha_4$. ($\alpha_1 + \alpha_2 + \alpha_3 + \alpha_4 = 2\pi$) and four folding angles $\rho_1, \rho_2, \rho_3, \rho_4$ whose corresponding inner creases are $c_1, c_2, c_3, c_4$. Figure 1 shows a developable degree-4 single-vertex creased paper and its two non-trivial rigidly folded states.

**Proposition 1** The configuration space of a developable degree-4 single-vertex creased paper is one of the following:

(a) If only one pair of adjacent inner creases are collinear (without loss of generality, suppose $c_1$ and $c_2$ are collinear), then $\rho_1 = \rho_2, \rho_3 = \rho_4 = 0$.

(b) If only one pair of opposite inner creases are collinear (without loss of generality, suppose $c_1$ and $c_3$ are collinear), the configuration space is the union of a line segment $\rho_1 = \rho_3, \rho_2 = \rho_4 = 0$ and a strictly monotonous 0-connected smooth curve in $\mathbb{R}^4$. The line segment and curve are symmetric to **0** and only intersect at **0**. The equation of both of these two branches can be derived by substituting $\alpha_1 + \alpha_4 = \pi$ and $\alpha_2 + \alpha_3 = \pi$ into equations (2) and (3) given below. Half of the second branch for $\rho_1 \geq 0$ is thus (at most differ by $2\pi$)

$$\rho_2 = 2 \arccos\left(\frac{\cos \alpha_2 \cos \xi - \cos \alpha_1}{\sin \alpha_2 \sin \xi}\right)$$
$$\rho_3 = -\rho_1 \quad (1)$$
$$\rho_4 = 2 \arccos\left(\frac{\cos \alpha_1 \cos \xi - \cos \alpha_2}{\sin \alpha_1 \sin \xi}\right)$$

where,
$$\xi = \arccos(\cos \alpha_1 \cos \alpha_2 - \sin \alpha_1 \sin \alpha_2 \cos \rho_1)$$

The domain of $\rho_1$ is a closed interval in $\mathbb{R}$. We call this a *straight-line* degree-4 single-vertex creased paper.

(c) If two pairs of inner creases are each collinear, the configuration space is the union of two line segments $\rho_1 = \rho_3, \rho_2 = \rho_4 = 0$ and $\rho_2 = \rho_4, \rho_1 = \rho_3 = 0$.

(d) If no pair of inner creases is collinear with every $\alpha_i \in (0, \pi)$, the configuration space is the union of two strictly monotonous 0-connected smooth curves in $\mathbb{R}^4$, which are symmetric to **0** and only intersect at **0**, whose expressions are:
   a) Half of branch 1 ($\rho_1 \geq 0$)

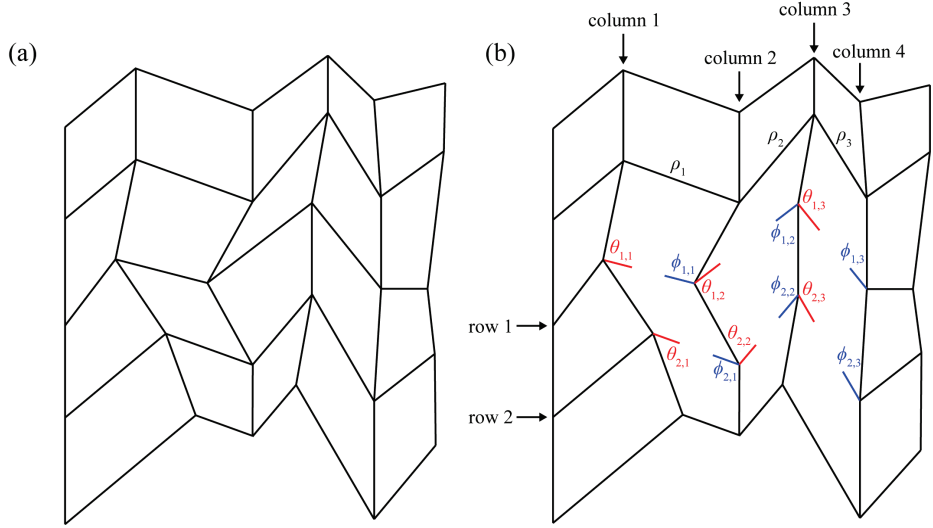

Figure 2: (a) shows the trivial rigidly folded state of a quadrilateral blanket, and (b) is its corresponding tree structure, where we "cut" inner creases connecting adjacent columns to make the interior of the crease pattern have no cycle. $\rho_j$ are the folding angles of the top row. $\theta_{i,j}$, $\phi_{i,j}$ and their corresponding inner creases are colored red and blue, respectively.

$$\rho_2 = \arccos\left(\frac{\cos\alpha_2 \cos\xi - \cos\alpha_1}{\sin\alpha_2 \sin\xi}\right) - \arccos\left(\frac{\cos\alpha_4 - \cos\alpha_3 \cos\xi}{\sin\alpha_3 \sin\xi}\right)$$
$$\rho_3 = \arccos\left(\frac{\cos\alpha_3 \cos\alpha_4 - \cos\xi}{\sin\alpha_3 \sin\alpha_4}\right) \quad (2)$$
$$\rho_4 = \arccos\left(\frac{\cos\alpha_1 \cos\xi - \cos\alpha_2}{\sin\alpha_1 \sin\xi}\right) - \arccos\left(\frac{\cos\alpha_3 - \cos\alpha_4 \cos\xi}{\sin\alpha_4 \sin\xi}\right)$$

b)  Half of branch 2 ($\rho_1 \geq 0$, at most differ by $2\pi$)
$$\rho_2 = \arccos\left(\frac{\cos\alpha_2 \cos\xi - \cos\alpha_1}{\sin\alpha_2 \sin\xi}\right) + \arccos\left(\frac{\cos\alpha_4 - \cos\alpha_3 \cos\xi}{\sin\alpha_3 \sin\xi}\right)$$
$$\rho_3 = -\arccos\left(\frac{\cos\alpha_3 \cos\alpha_4 - \cos\xi}{\sin\alpha_3 \sin\alpha_4}\right) \quad (3)$$
$$\rho_4 = \arccos\left(\frac{\cos\alpha_1 \cos\xi - \cos\alpha_2}{\sin\alpha_1 \sin\xi}\right) + \arccos\left(\frac{\cos\alpha_3 - \cos\alpha_4 \cos\xi}{\sin\alpha_4 \sin\xi}\right)$$

where,
$$\xi = \arccos(\cos\alpha_1 \cos\alpha_2 - \sin\alpha_1 \sin\alpha_2 \cos\rho_1)$$
The domain of $\rho_1$ is a closed interval in $\mathbb{R}$.

(e) If the sector angles are not in the cases (a) − (d), the configuration space is trivial (Abel et al. [4]).

*Proof.* We can first assume $\alpha_1 + \alpha_2 \geq \alpha_3 + \alpha_4$, equations (2) and (3) can be derived directly from Figure 1 with the spherical cosine theorem, then by verification we know they are also satisfied when $\alpha_1 + \alpha_2 < \alpha_3 + \alpha_4$. In that case the value of folding angles will differ by $2\pi$ in equation (3). Note that $\xi \in (0, \pi)$ because every $\alpha_l \in (0, \pi)$. We need to show that the smooth curves in case (d) and the second branch of case (b) are strictly monotonous, which can be verified by direct calculation.

## 3. General Theory

For a developable quadrilateral creased paper, we can always find the smallest creased paper covering single-vertex creased papers incident to each group of connected inner panels, called a quadrilateral blanket. Generically, only quadrilateral blankets affect the rigid-foldability [1, 2], therefore we just consider a quadrilateral blanket. Without loss of generality, we suppose a quadrilateral blanket is generated by a row and several columns of developable degree-4 single-vertex creased papers. The following theorem gives a sufficient and necessary condition for the trivial rigidly folded state of a quadrilateral blanket to be rigid-foldable.

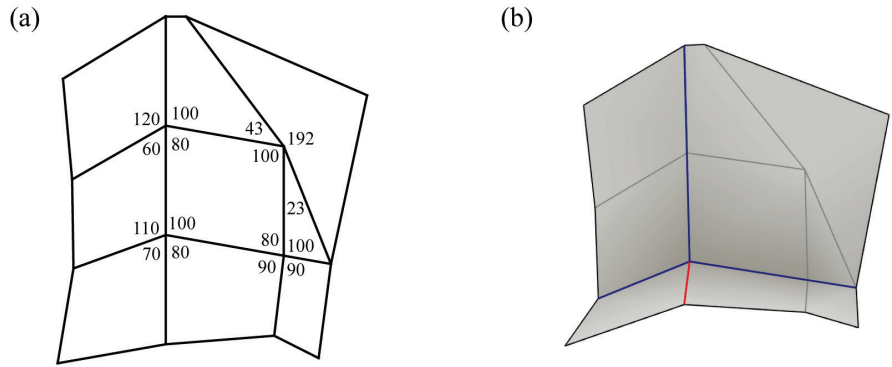

Figure 3: (a) A rigid-foldable creased paper, where not all its single-vertex creased papers are rigid-foldable. The magnitude of sector angles are labeled in degrees. (b) One of its rigidly folded state, plotted by Freeform Origami (Tachi [6]). The mountain and valley creases are colored red and blue, while the grey creases are not folded.

**Theorem 1** [2] Consider the trivial rigidly folded state of a quadrilateral blanket. If and only if

1) the corresponding rigidly folded state of the tree structure (see Figure 2(b)) is rigid-foldable in a closed interval of all the folding angles.
2) for all $i,j$, the following equation is satisfied simultaneously in this closed interval.

$$\theta_{i,j} \equiv \phi_{i,j} \tag{4}$$

then this trivial rigidly folded state is rigid-foldable in this closed interval. Figure 2 defines $\theta_{i,j}$ and $\phi_{i,j}$.

**Remark 1** In a quadrilateral rigid-foldable creased paper, it is not necessary for all degree-4 single-vertex creased papers to be rigid-foldable. There may be a non-rigid-foldable single-vertex creased paper whose panels do not have relative rigid folding motions. An example is shown in Figure 3, which also explains why Theorem 1 is wider than the condition described in Tachi [3]. Besides, unlike the flat-foldable case, a non-trivial rigidly folded state cannot guarantee a rigid folding motion for a general quadrilateral creased paper.

**Corollary 1.1** [2] To satisfy equation (4), one sufficient condition is, for all $i,j$,

$$\theta_{i,j} \equiv \phi_{i,j} = \pm\rho_j \tag{5}$$

and one way to realize it is to use the tessellation of "units".

## 4. Units

**Definition 2.** As shown in Figure 4, a *unit* is the union of two degree-4 single-vertex creased papers sharing two panels, and the folding angles on the same side are constantly of the same magnitude in a closed interval, i.e.

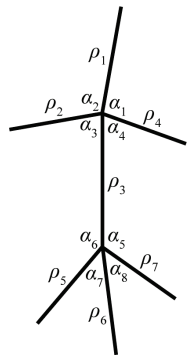

Figure 4: A unit with its 8 sector and 7 folding angles. Here we only draw the inner creases.

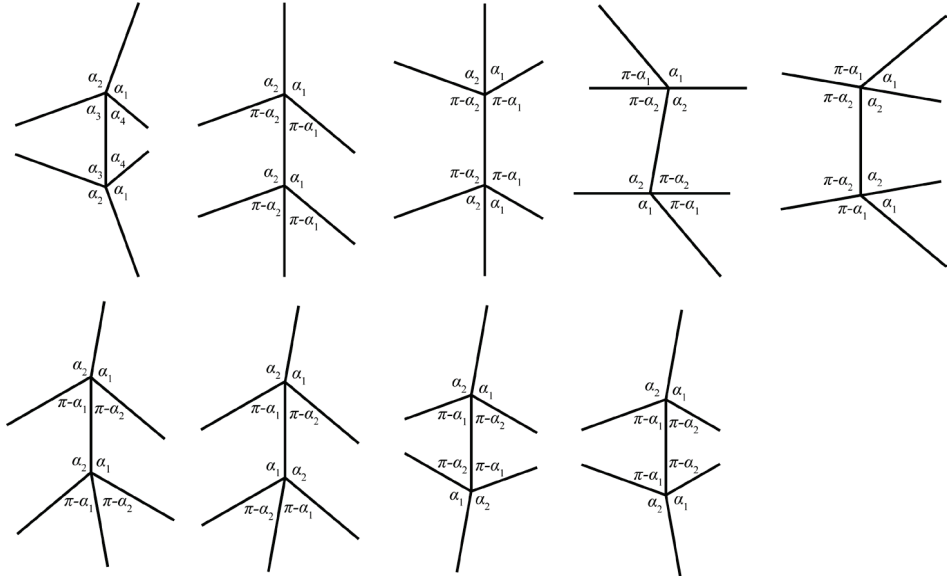

Figure 5: The nine basic units that are discussed in detail in [2], $\alpha_i$ are independent sector angles.

$$\rho_2 \equiv \pm\rho_5, \quad \rho_4 \equiv \pm\rho_7 \tag{6}$$

**Remark 2.** We also introduce nine "basic units" in [2], as shown in Figure 5, which satisfy equation (6) and contain two identical vertices. Although using identical vertices is not necessary, it makes the relation of sector angles in basic units simpler, and reduces the complexity in design.

Without loss of generality, for the top vertex, we label the sector angles adjacent to the inner crease connecting these two vertices by $\alpha_3$ and $\alpha_4$, then label $\alpha_1$ and $\alpha_2$ counter-clockwise. Here we do not consider cases (a), the first branch of (b), (c) and (e) of Proposition 1, because a single-vertex creased paper in these cases either has the same rigid folding motion as two panels incident to an inner crease or is not rigid-foldable.

From equations (2) and (3), we know generally the relation among the sector and folding angles of a degree-4 single-vertex creased paper is coupled, which makes it hard to find special solutions for equation (6). Here we only discuss a simpler case, where a unit consists of two flat-foldable vertices. Note that by choosing two straight-line vertices in Proposition 1 we can only get basic units.

**New Flat-foldable Unit**

Here we introduce a new flat-foldable unit to go alongside the basic unit previously defined. For a flat-foldable degree-4 single-vertex creased paper, where $\alpha_1 + \alpha_3 = \pi$, $\alpha_2 + \alpha_4 = \pi$, $\alpha_1, \alpha_2 \in (0, \pi)$ and not simultaneously equal to $\pi/2$. Equation (2) and (3) degenerate to (Tachi and Hull [5]):

a) Branch 1

$$\tan\frac{\rho_2}{2} = \left(\frac{\sin\frac{\alpha_2 - \alpha_1}{2}}{\sin\frac{\alpha_2 + \alpha_1}{2}}\right) \tan\frac{\rho_1}{2}$$
$$\rho_3 = \rho_1$$
$$\rho_4 = -\rho_2 \tag{7}$$

b) Branch 2

$$\tan\frac{\rho_2}{2} = -\left(\frac{\cos\frac{\alpha_2 - \alpha_1}{2}}{\cos\frac{\alpha_2 + \alpha_1}{2}}\right) \tan\frac{\rho_1}{2}$$
$$\rho_3 = -\rho_1$$
$$\rho_4 = \rho_2 \tag{8}$$

As shown in Figure 6, if we want a flat-foldable unit to satisfy equation (6), there are four modes. Equation (6) can be written as:

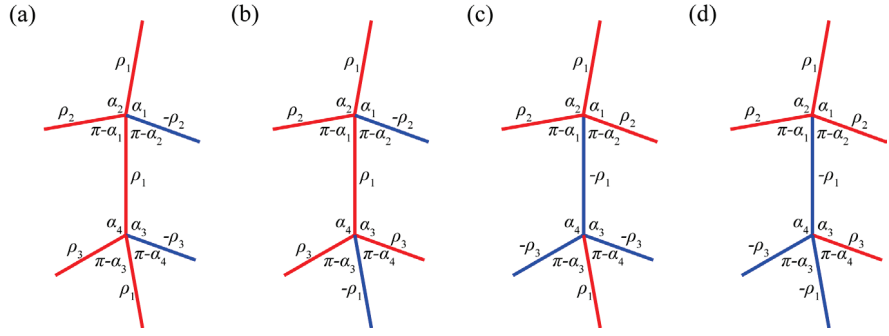

Figure 6: (a) − (d) show four branches of rigid folding motion of the new flat-foldable unit introduced in this article. We color the mountain creases red, while the valley creases blue. Note that swapping the mountain-valley assignment results in symmetric rigid folding motions in the same branch. We show that, in fact (b) and (d) are not possible.

$$\frac{\tan\frac{\alpha_2}{2} - \tan\frac{\alpha_1}{2}}{\tan\frac{\alpha_2}{2} + \tan\frac{\alpha_1}{2}} = \pm \frac{\tan\frac{\alpha_4}{2} - \tan\frac{\alpha_3}{2}}{\tan\frac{\alpha_4}{2} + \tan\frac{\alpha_3}{2}} \tag{9a}$$

$$\frac{\tan\frac{\alpha_2}{2} - \tan\frac{\alpha_1}{2}}{\tan\frac{\alpha_2}{2} + \tan\frac{\alpha_1}{2}} = \pm \frac{1 + \tan\frac{\alpha_4}{2}\tan\frac{\alpha_3}{2}}{1 - \tan\frac{\alpha_4}{2}\tan\frac{\alpha_3}{2}} \tag{9b}$$

$$\frac{1 + \tan\frac{\alpha_2}{2}\tan\frac{\alpha_1}{2}}{1 - \tan\frac{\alpha_2}{2}\tan\frac{\alpha_1}{2}} = \pm \frac{1 + \tan\frac{\alpha_4}{2}\tan\frac{\alpha_3}{2}}{1 - \tan\frac{\alpha_4}{2}\tan\frac{\alpha_3}{2}} \tag{9c}$$

$$\frac{1 + \tan\frac{\alpha_2}{2}\tan\frac{\alpha_1}{2}}{1 - \tan\frac{\alpha_2}{2}\tan\frac{\alpha_1}{2}} = \pm \frac{\tan\frac{\alpha_4}{2} - \tan\frac{\alpha_3}{2}}{\tan\frac{\alpha_4}{2} + \tan\frac{\alpha_3}{2}} \tag{9d}$$

Each equation in (9a) − (9d) represents one branch of rigid folding motion that keeps equation (6) satisfied. Because here every $\alpha_i \in (0, \pi)$, $\tan\frac{\alpha_i}{2} \in (0, +\infty)$. Moreover, functions $y = \frac{1+x}{1-x}$ and $y = \frac{1-x}{1+x}$ when $x > 0$ are (piecewise) monotonous and have no intersection, and their ranges are symmetric to 0, thus we know equations (9b) and (9d) have no solution. Equations (9a) and (9c) can be simplified as:

$$\frac{\tan\frac{\alpha_1}{2}}{\tan\frac{\alpha_2}{2}} = \frac{\tan\frac{\alpha_3}{2}}{\tan\frac{\alpha_4}{2}} \quad \text{or} \quad \frac{\tan\frac{\alpha_1}{2}}{\tan\frac{\alpha_2}{2}} = \frac{\tan\frac{\alpha_4}{2}}{\tan\frac{\alpha_3}{2}} \tag{10a}$$

$$\tan\frac{\alpha_1}{2}\tan\frac{\alpha_2}{2} = \tan\frac{\alpha_3}{2}\tan\frac{\alpha_4}{2} \quad \text{or} \quad \tan\frac{\alpha_1}{2}\tan\frac{\alpha_2}{2}\tan\frac{\alpha_3}{2}\tan\frac{\alpha_4}{2} = 1 \tag{10b}$$

Therefore the branches of rigid folding motion shown in Figure 6(b) and 6(d) cannot satisfy equation (6), and in Figure 6(a) and 6(c), one set of sector angles $\alpha_1, \alpha_2, \alpha_3$ can always determine two possible $\alpha_4$. There are three independent sector angles in this flat-foldable unit.

## 5. Some Examples

In this section we continue to use the idea of "stitching" first described in [2]. Based on the known nine basic units in [2] and the flat-foldable units introduced in this article, we can generate a family of developable quadrilateral creased papers that are rigid-foldable and have only one degree of freedom during the rigid folding motion, except at the trivial rigidly folded state. Examples are shown in Figures 7 and 8.

When designing a crease pattern, we naturally focus on the degrees of freedom we have in choosing sector angles. Every new flat-foldable unit has three, while stitching in a column will increase the number by one, and we can calculate the total independent sector angles unit by unit, unlike the case where only basic units are involved. This is because stitching in a column can add more independent sector angles. In transverse direction, if there is a row of non-parallel creases, the amount of independent sector angles decrease by the number of inner panels of this row, otherwise for a row of parallel creases, the amount will not decrease. This is because the sum of sector angles of all the inner panels should be $2\pi$. Note that the number of independent sector angles should be greater than or equal to zero when designing the whole crease pattern. In general, less independent sector angles will increase the difficulty in design because the sector angles will be coupled by more constraints.

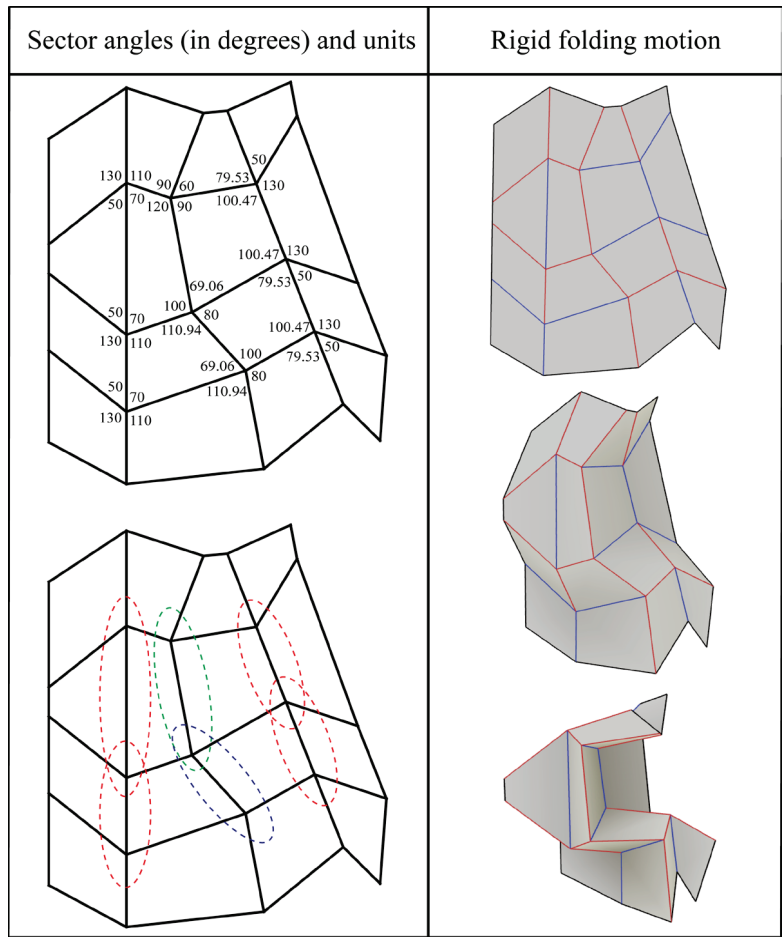

Figure 7: A developable quadrilateral creased paper that contains four straight-line basic units (in red), a flat-foldable basic unit (in blue) and a flat-foldable unit (in green). We show the sector angles and where the units locate separately. The rigid folding motions are plotted by Freeform Origami (Tachi [6]) from the trivial state to the final state. The mountain and valley creases are colored red and blue. Independent sector angles: $2 + 3 + 2 - 2 = 5$. Branches of rigid folding motion: $1 \times 1 \times 1 = 1$.

For a developable quadrilateral creased paper, there may be more than one branch of rigid folding motion, corresponding with different mountain-valley assignments. Unlike a flat-foldable basic unit, every flat-foldable unit obtained by solving equation (10a) or (10b) has only one branch of rigid folding motion. Then we can calculate the number of branches of rigid folding motion by multiplying the number of each column.

## 6. Conclusion

In this paper, we have developed flat-foldable units that can be stitched with the known basic units to construct new 1-DOF rigid-foldable developable quadrilateral creased papers. Compared to basic units, these flat-foldable units have more complex relations among their sector angles, which can generate more irregular rigid folding motions. These analyses can also be used to generate new non-developable quadrilateral creased papers mentioned in [2].

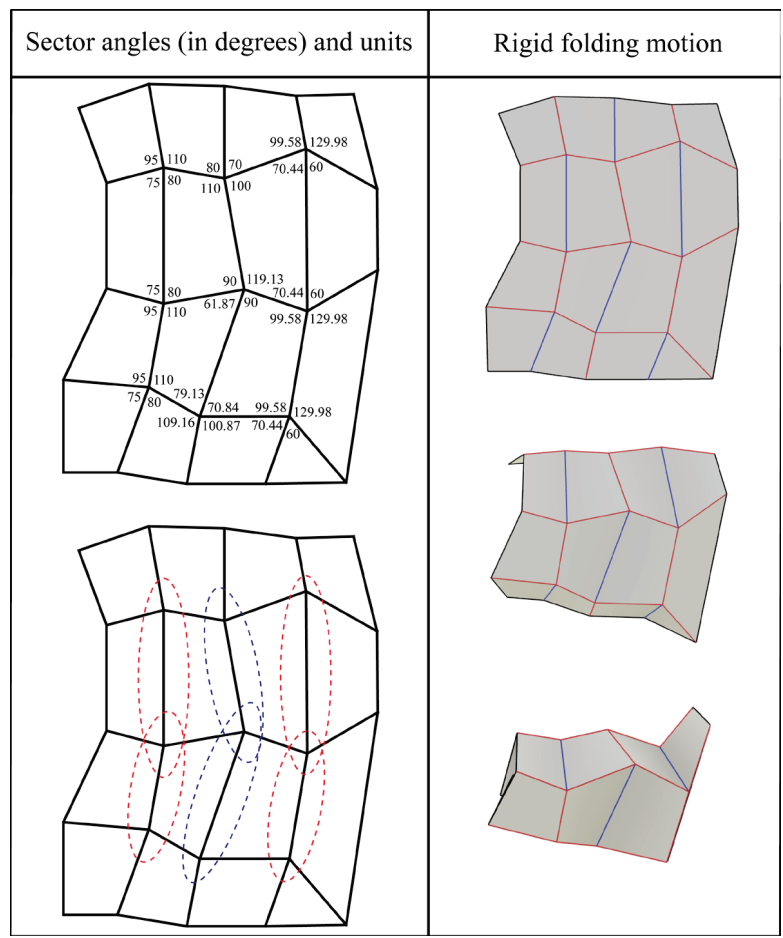

Figure 8: A developable quadrilateral creased paper that contains four general basic units (in red) and two flat-foldable units (in blue). We show the sector angles and where the units locate separately. The rigid folding motions are plotted by Freeform Origami (Tachi [6]) from the trivial state to the final state. The mountain and valley creases are colored red and blue. Independent sector angles: $3 + 3 + 1 + 3 - 4 = 6$. Branches of rigid folding motion: $2 \times 1 \times 2 = 4$.